\documentclass[11pt]{amsart}

\usepackage{hyperref}

\makeatletter
\def\@normalsize{\@setsize\normalsize{10pt}\xpt\@xpt
\abovedlayskip 10pt plus2pt minus5pt\belowdisplayskip
\abovedisplayskip \abovedisplayshortskip \z@
plus3pt\belowdisplayshortskip 6pt plus3pt
minus3pt\let\@listi\@listI}
\def\subsize{\@setsize\subsize{12pt}\xipt\@xipt}
\def\section{\@startsection {section}{1}{\z@}{1.0ex plus 1ex minus
2ex}{.2ex plus .2ex}{\large\bf}}
\def\subsection{\@startsection {subsection}{2}{\z@}{.2ex plus 1ex}
{.2ex plus .2ex}{\subsize\bf}} \makeatother

\newtheorem{lemma}{Lemma}

\newtheorem{theorem}{Theorem}

\newtheorem{question}{Question}

\theoremstyle{remark}

\begin{document}

\baselineskip=12pt

\email{alexei.ostrovski@gmx.de}
\address{Follow me on  \href{https://www.researchgate.net/profile/Alexey_Ostrovsky/}{ResearchGate} }



\title{\bf { Open--constructible  functions \\(Corrected version, January 2014)
}}

\author {\href{https://www.researchgate.net/profile/Alexey_Ostrovsky/}{Alexey Ostrovsky}}
 
\maketitle

\begin{abstract}     Let $f$ be a continuous  function between  subspaces $X,Y$ of   the Cantor set  $ \textbf{C}.$
We prove that:

 if $f$ is  one-to-one  and maps
 open sets into resolvable,  then $f$ is a piecewise homeomorphism  and 

  if  $f$   maps discrete subsets  into  resolvable, 
    then $f$ is  piecewise open.  
\end{abstract}

   \vspace{0.15in }

 \section{Introduction}
 
 \vspace{0.15in }

  The present paper continues the series of publications 
about decomposibility of Borel functions \cite{CC}, \cite{O3} -  see also \cite{GM}, \cite{HO}  where functions of such type are the main subject.

   A subset  $E$ of a topological space $X$  is   \emph{resolvable} if  for each  nonempty closed in $X$ subset   $F$  we have
   
   $$ cl_X(F \cap E )  \cap  cl_X(F \setminus E ) \not = F  $$

   \vspace{0.15in }

Recall that a function $f$  is   open    if  it
maps   open  sets into open ones.   
More generally, a function $f$  is said to be \emph{open--resolvable}  (resolvable in \cite{SG}) if   $f$ maps   open  sets into resolvable ones.

     \vspace{0.10in }
 A function $f : X \to  Y$  for which   $X \subset  \textbf{C}$ admits a  countable,  closed and disjoint  cover $ \mathcal{C}$,   such that for each   $C    \in  \mathcal{C} $ the 
restriction $f | C$ is open,  is  called    \emph{piecewise open}.

    \vspace{0.10in }

\begin{theorem}  

  Let $f$  be  a  continuous,   one-to-one and   open--resolvable function   between    $X,Y \subset  \textbf{C}.$ Then $f$ is  piecewise open and hence  is piecewise homeomorphism.

  \end{theorem}
    
 \vspace{0.05in }      
    



Note, that using  standard sets $H,$ we obtain  the proof  for   open-resolvable functions     simpler  than  for  open-constructible. 

 Analogously,  using   $H$-sets we can easy  extend the proof  for   closed-constructible functions   in \cite{CC} to the case of closed-resolvable.

      \vspace{0.10in }  
     \noindent
     Standard set $H$  will be  used in the proof of the following theorem:

   
    \begin{theorem}  

  If    a  continuous    function  $f: X \to Y $  between    $X,Y \subset  \textbf{C}$  maps discrete subsets in $X$ onto  resolvable, 
    then $f$ is  piecewise open.    

  \end{theorem}

    \vspace{0.25in }

   \vspace{0.15in }
   
\subsection {Standard set $H$.}  

 \vspace{0.15in }  

The  set $H$ was introduced by W. Hurewicz   \cite{H} and has a lot of applications. 
 $H$ (called standard in \cite{AVO},\cite{A})  is a countable set without isolated points:
 
  $$H =  \{ p,...., p_{i_1,...i_k}: k\in N^+;  i_1,...,i_k     \in N^+  \} $$       such that 
  $$p_{i_1}  \to  p,  \text{ as }        i_{i_1} \to   \infty$$

   $$p_{i_1,...i_k, i_{k+1}}  \to  p_{i_1,...i_k},   \text{ as } i_{k+1}  \to   \infty$$

    \vspace{0.10in }

 Obviously, $H$ is homeomorphic to the space of rational  $ \textbf{Q}.$

 Using the metric in $X  \supset  H$ we can suppose  additionally 
  that there are  decreasing  bases $U^i(p),$ and $U^i(p_{i_1,...i_k}) $  
    at points $p$ and      $p_{i_1,...i_k} $   satisfying 
      conditions a), b)  and c) below:
      
      \vspace{0.10in } 
      
 a)  $U^{i_k}(p_{i_1,...i_k})  \supset U^{1}î(p_{i_1,...i_{k+1}}),  i_k   \in N^+$
 
    \vspace{0.10in }

 b)  for   $i_k^{ \prime }  \not = i_k^{ \prime \prime }$ we have $U^1î(p_{i_1,...i_k^{ \prime }}) \cap U^{{1}}î(p_{i_1,...i_{i_k^{ \prime\prime }}}) =\emptyset . $
 
     \vspace{0.10in } 

  c)  diam($U^{i_1}(p_{i_1,...i_k}) )< 1/({i_1+...+i_k)}$  
  
  
   \vspace{0.15in }

\section{ Construction of  $Z \subset X$   for which $f|Z$ is nowhere open}

       \vspace{0.15in }

        Given  a function   $f : X  \to Y$, we shall construct  in the next Lemma 1  a subset $Z  \subset X$     on which the restriction $f|Z$ is   \emph{nowhere open} on $Z$; i.e. 
 for every  clopen in $X$   subset $U$  
     the restriction  $f|(U \cap Z)$ is    not open.

 \vspace{0.05in }

      \begin{lemma}  Let  $f : X  \to Y$  be a continuous  function from a subspace   $X$ of the Cantor set  $ \textbf{C}$ onto a   metrizable space  $Y.$  
 Then  there is a closed subset   $Z  \subset X $  such that 
the restriction $f |Z$ is nowhere open  on $Z$ and the restriction  $f|(X  \setminus Z)$  is piecewise open.

     \vspace{0.15in }    
 \end{lemma}

\noindent
Proof of Lemma 1.  Let us begin by proving the first part of the assertion from lemma  stating that for some $Z$ the restriction $f |Z$ is nowhere open  on $Z.$  Indeed, if  for some  nonempty clopen set $V
\subset X $ the restriction  $f|V$ is open, then  we could construct the closed set
$$X_1 = X   \setminus V $$  and the corresponding
restriction
 $$f|X_1: X_1 \to f(X_1) .$$

\noindent
Repeating  this process, we could also construct a chain of
closed sets ($X_{\gamma} = \bigcap_{\beta < \gamma} X^{\beta}$ for a limit $\gamma$)

$$X \supset X_1 \supset ...\supset X_{\gamma} \supset...$$
which, as we know, stabilizes
  at some   $\gamma_0  < \omega_1$. Therefore, there exists a subspace $Z$ for which holds true

  $$Z =  X_{\gamma_0} =  X_{ \gamma_0+1}= ...$$
and  the restriction $f|Z$ is   nowhere open  on $Z.$

The second part of Lemma 1 stating that $f|(X  \setminus Z)$   is piecewise open, obviously, satisfied.

  In what follows, it is convenient in the case when $Z$  is empty to regard a    piecewise open  function   $f|X =f|(X  \setminus Z)$   onto $Y.$
           
\qed

   \vspace{0.10in }

 \subsection{Proof of Theorems 1 and 2. }

    \vspace{0.10in }

 On the step 1   take  a point $x \in   Z$   
and
  a  base  of clopen  in $X$ neighborhoods   $U^k(x)   \subset  X$ with diametr less than $1/k.$
  
  Since $f$ is  nowhere open on $Z,$
      there  are      $$U^l(x_{ k})  \subset     Z   \cap (U^k(x)   \setminus  (U^{k+1}(x) )$$  
     

    such that   $x_{ k}    \to  x,$ $f(x_{ k}) = p_k    \to p = f(x)$   and

     $$ f^{-1}(p_{ k,l} )   \cap  U^1(x_{ k})    =   \emptyset $$

       
       Take    $  x_{ k,l}   \in U(x_{ k,l})   \subset       \cap (U^l(x_k)   \setminus  (U^{l+1}(x_k) ).$    
       
     \vspace{0.10in }

     For the basic induction step $m$   we  suppose that the points $p = f(x), .... p_{i_1,...i_k} = f(x_{i_1,...i_k})$  satisfying the conditions a), b),c) of definition the set $H$  are constructed. 
     
     Analogously   
      pick in some  clopen sets  $U^1(x_{i_1,...i_k})  $ the points  $x_{ k}   \to x $ with  pairwise disjunct   clopen   neighborhoods  $U^1(x_{ k}) $  such  that  for some sequence 
        $p_{ k,l}    \to  p_{ k} = f(x_k)$   we have

       $$ f^{-1}(p_{ k,l} )   \cap  U^1(x_{ k})    =   \emptyset $$       
       
   \noindent     
  Take the points  $$x_{ k,l}     \in     f^{-1}(p_{ k,l} )   \setminus  U^1(x_k)  $$  
     and  clopen  sets $$U^1(x_{k,l} )  \subset X  \setminus  U^1(x_k).   $$ 
    
    Since $f$ is continuous we can suppose that  $U^1(x_{k,l} ) $  are disjunct with all  $U^1(x_{k}). $ 

Since $f$ is   nowhere open on   $ Z   \cap U^1(x_{ k,l})  $ we  can  repeat the construction for $U^1(x_{k,l} )$  etc.
Analogously   we obtain


    $$p_{i_1,...i_m, k,l}  \to  p_{i_1,...i_m, k} \text{ as } l  \to   \infty  $$

     $$x_{i_1,...i_m, k,l} \in   f^{-1} (p_{i_1,...i_m, k,l})  \setminus  U^1(x_{i_1,...i_m, k}). $$

   $$U^1(x_{i_1,...i_m, k}) \cap f^{-1} (p_{i_1,...i_m, k,l})  = \emptyset,  l   \in N^+$$

We can suppose that $U^1(p_{i_1,...i_m, k})$ are  clopen and  there is  a clopen neighborhood   $U^1(x_{i_1,...i_m, k,l})$  of point $x_{i_1,...i_m, k,l}, $ disjunct with all $U^1(x_{i_1,...i_m, k})$.

     
Denote       $$D =   \{ x_{i_1}, x_{i_1,i_2,i_3}\dots x_{i_1, \dots i_{2k+1} }: k \in N^+;  i_1,...,i_{2k+1}    \in N^+  \} $$    
 

  By our construction $D$ is discrete and $f(D)$  is dense and codense in $$H =   \{  p_{i_1, \dots i_{k} }: k \in N^+;  i_1,...,i_{k}    \in N^+  \} $$ 
  
  \noindent
  that  proves Theorem 2.
  
  To prove Theorem  1  we note, that in case of one-to-one functions 
   the open set  $O(D) =  \{ U^1(x_{i_1}), \dots U^1(x_{i_1, \dots i_{2k+1}) }: k \in N^+;  i_1,...,i_{2k+1}    \in N^+  \}   $ has the same image as  $D$.
   
     \qed   
   \vspace{0.15in }
        
       \begin{question}Are the continuous open--$\Delta^0_2$ functions       
        (even for Polish  or analytic  spaces $X \subset   \textbf{C}$) piecewise or countably open?
         \end{question} 
         
 \vspace{0.15in }
 
     \vspace{0.20in }

 \vspace{0.30in }

\end{document}